\documentclass[12pt]{article} 
\usepackage[all,cmtip]{xy}
\usepackage{amsmath,amssymb,eucal,bbm} 
\font\tenrm=cmr10

\font\cmssl=cmss10 at 12 pt  
   
\font\bigss=cmssdc10 scaled 2300

\font\cmsslll=cmss10 at 14 pt

\renewcommand{\a}{\alpha}  
\renewcommand{\b}{\beta}  
  
\renewcommand{\d}{\delta}

\newcommand{\g}{\gamma}

\renewcommand{\o}{\omega}  
  
\renewcommand{\r}{\rho}

\newcommand{\G}{\Gamma}

\renewcommand{\O}{\Omega}

\newcommand{\bC}{\mathbb{C}}  
\newcommand{\bR}{\mathbb{R}}  
  
\newcommand{\bH}{\mathbb{H}}

\newcommand\SO{\mathrm{SO}}

\newcommand{\id}{{\mathbbm{1}}}   
    
\newcommand{\p}{\partial}  
    
\renewcommand{\square}{\kern1pt\vbox  
               {\hrule height 0.6pt\hbox{\vrule width 0.6pt\hskip 3pt  
    \vbox{\vskip 6pt}\hskip 3pt\vrule width 0.6pt}\hrule height0.6pt}  
    \kern1pt}  

\newcommand{\ra}{\rightarrow}

\newcommand{\wt}{\widetilde}

\newtheorem{Pb}{Problem} 
   
\newtheorem{Ex}{Example}   
\newtheorem{Th}{Theorem}  
\newtheorem{Prop}{Proposition}  
  
\newtheorem{Cor}{Corollary}  
\newtheorem{Lem}{Lemma}  
\newtheorem{Def}{Definition} 
\newcommand{\bP}{\begin{Pb}\ \ } 
\newcommand{\eP}{\end{Pb}}  
\newcommand{\bt}{\begin{Th}\ \ }  
\newcommand{\et}{\end{Th}}  
\newcommand{\bp}{\begin{Prop}\ \ }  
\newcommand{\ep}{\end{Prop}}  
\newcommand{\bc}{\begin{Cor}\ \ }  
\newcommand{\ec}{\end{Cor}}  
\newcommand{\bl}{\begin{Lem}\ \ }  
\newcommand{\el}{\end{Lem}}  
\newcommand{\bd}{\begin{Def}\ \ }  
\newcommand{\ed}{\end{Def}}  
  
\newcommand{\pf}{\noindent{\it Proof:\ \ }}  
\newcommand{\qed}{\hfill\square}  
\newcommand{\n}{\nabla}

\newcommand{\be}{\begin{equation}}  
\newcommand{\ee}{\end{equation}}  
  
\newcommand\re[1]{(\ref{#1})}  
\newcommand{\arr}{\begin{array}{rlll}}  
\newcommand{\ea}{\end{array}}  
\newcommand{\bea}{\begin{eqnarray}}  
\newcommand{\eea}{\end{eqnarray}}  
\newcommand{\bean}{\begin{eqnarray*}}  
\newcommand{\eean}{\end{eqnarray*}}  
  
\catcode`@=11  
\@addtoreset{equation}{section}  
\catcode`@=12

\begin{document}  
 \rightline{} 
\vskip 1.5 true cm  
\begin{center}  
{\bigss  Conification of K\"ahler and hyper-K\"ahler manifolds}\\[.5em]
\vskip 1.0 true cm   
{\cmsslll  D.V.\ Alekseevsky$^1$, V.\ Cort\'es$^2$ and T.\ Mohaupt$^3$} \\[3pt] 
$^1${\tenrm Institute for Information Transmission Problems\\
B.\ Karetny per.\ 19\\ 
101447 Moscow, Russia\\
daleksee@staffmail.ed.ac.uk}\\[1em]  
$^2${\tenrm   Department of Mathematics\\  
and Center for Mathematical Physics\\ 
University of Hamburg\\ 
Bundesstra{\ss}e 55, 
D-20146 Hamburg, Germany\\  
cortes@math.uni-hamburg.de}\\[1em]  
$^3${\tenrm Department of Mathematical Sciences\\ 
University of Liverpool\\
Liverpool L69 3BX, UK\\  
Thomas.Mohaupt@liv.ac.uk}\\[1em] 
July 13, 2012 
\end{center}  
\vskip 1.0 true cm  
\baselineskip=18pt  
\begin{abstract}  
\noindent  
Given a K\"ahler manifold $M$ endowed with a 
Hamiltonian Killing vector field $Z$,  
we construct a conical K\"ahler
manifold $\hat{M}$ such that $M$ is recovered as a K\"ahler
quotient of $\hat{M}$. Similarly, given a hyper-K\"ahler manifold
$(M,g,J_1,J_2,J_3)$ endowed with a 
Killing vector field $Z$, Hamiltonian with respect to the 
K\"ahler form of $J_1$ and satisfying 
$\mathcal{L}_ZJ_2= -2J_3$, we construct a hyper-K\"ahler cone $\hat{M}$
such that $M$ is a certain hyper-K\"ahler quotient of $\hat{M}$. 
In this way, we recover a theorem by Haydys. 
Our work is motivated by the problem of relating the supergravity
c-map to the rigid c-map. We show that any hyper-K\"ahler manifold
in the image of the c-map admits a Killing vector field 
with the above properties. Therefore, it gives rise to a hyper-K\"ahler 
cone, which in turn defines a quaternionic K\"ahler manifold. 
Our results for the signature of the metric and the sign of the 
scalar curvature are consistent with what we know about the supergravity
c-map. 
\end{abstract}

\section*{Introduction}
Let us recall that there is an interesting geometric construction
called the c-map, which was found by theoretical physicists. 
There are in fact two versions of the c-map: the 
supergravity c-map and the rigid c-map.  
The \emph{supergravity c-map} associates a quaternionic K\"ahler
manifold of negative scalar curvature 
with any projective special K\"ahler manifold, see 
\cite{FS,H2,CM}. The metric is explicit but rather 
complicated. The \emph{rigid c-map} is much simpler, see
\cite{CFG,H1,ACD}. It associates a 
hyper-K\"ahler manifold with any affine 
special K\"ahler manifold. 
The initial motivation for this work was our idea to reduce 
the supergravity c-map to the rigid c-map
by means of a conification of the hyper-K\"ahler 
manifold obtained from the rigid c-map. Let us explain this idea
in more detail. 

Since any projective special K\"ahler manifold
$\bar{M}$ is the base of a $\bC^*$-bundle with the total space
a conical affine special K\"ahler manifold $M$, we have the
following diagram:
\[
\xymatrix{
M^{2n} \ar@{|->}[r]^c \ar@{->}[d]_{\bC^*}& N^{4n} \\
\bar{M}^{2n-2} \ar@{|->}[r]^{\bar{c}} & \bar{N}^{4n} \\
}
\]
where $c$ stands for the rigid c-map, $\bar{c}$ for the 
supergravity c-map and $N$, $\bar{N}$ are the 
resulting (pseudo-)hyper-K\"ahler and quaternionic
K\"ahler manifolds, respectively. We have indicated the real 
dimension. This shows that $N$ cannot simply be the 
Swann bundle $\hat{N}$ over  $\bar{N}$. 
In fact, $N$ is in general 
not conical and the (pseudo-)hyper-K\"ahler cone 
$\hat{N}$ should be obtained from $N$
by a certain conification procedure $N^{4n}\stackrel{con}{\mapsto}
\hat{N}^{4n+4}$ such that the following diagramm commutes:
\[
\xymatrix{
M^{2n} \ar@{|->}[r]^c \ar@{->}[d]_{\bC^*} & N^{4n} \ar@{|->}[r]^-{con} &
\hat{N}^{4n+4} \ar@{->}[d]^{\bH^* / \pm 1} \\
\bar{M}^{2n-2} \ar@{|->}[rr]^{\bar{c}} & & \bar{N}^{4n} \\
}
\]
We are also interested in the analogous problem for the r-map, 
where we have a diagramm of the form:
\[
\xymatrix{
M^{n} \ar@{|->}[r]^r \ar@{->}[d]_{\bR^{>0}} & N^{2n} \ar@{|->}[r]^-{con} &
\hat{N}^{2n+2} \ar@{->}[d]^{\bC^*} \\
\bar{M}^{n-1} \ar@{|->}[rr]^{\bar{r}} & & \bar{N}^{2n} \\
}
\]
Now $M$ is an affine special real manifold with homogeneous 
cubic prepotential, $\bar{M}$ is the corresponding  projective
special real manifold, $r$ is the rigid r-map \cite{CMMS1,AC}, 
$\bar{r}$ is the supergravity r-map \cite{DV,CM} and
$\hat{N}$ is the conical affine special K\"ahler
manifold over the projective special K\"ahler manifold
$\bar{N}$. 

An important inspiration for our work has been the paper 
\cite{Haydys} by Haydys, see also \cite{APP} in which 
Haydys construction is called \emph{QK/HK correspondence}.
The construction has two parts. The first part is 
the hyper-K\"ahler reduction of a hyper-K\"ahler cone 
with respect to a Hamiltonian Killing vector field which is 
compatible with the cone structure. The  hyper-K\"ahler
manifold $(M,g,J_1,J_2,J_3)$ obtained by such a reduction inherits 
a Killing vector 
field $Z$ which preserves one of the three complex structures
$J_1$  of the hyper-K\"ahler triplet ($J_\a$) and rotates the two other ones. 
The second part is the inversion of the reduction, which is 
much more involved than the first part.  
As a result of our careful analysis, we are able to
give our own proof of the inversion recovering and extending the results 
by Haydys.  Under the assumptions stated precisely in Section \ref{HKSec}, 
the conical hyper-K\"ahler structure 
is rigorously established in Theorem \ref{hKconifThm}. The final formulas 
are explicit enough to allow for further progress 
in the study of hyper-K\"ahler manifolds obtained by 
such a conification. As an example, we can easily compute 
the signature and scalar curvature 
of the resulting quaternionic K\"ahler manifolds, see Corollary
\ref{signatureCor} and Corollary \ref{RiemCor}. These results 
are new even in the case when the initial hyper-K\"ahler metric
is positive definite, as considered in \cite{Haydys}. We show that
(positive definite) quaternionic K\"ahler manifolds of negative  
scalar curvature can be obtained from indefinite 
as well as from positive definite hyper-K\"ahler manifolds,
whereas quaternionic K\"ahler manifolds of positive   
scalar curvature do always require a positive definite 
initial metric.

We prove that a similar, but simpler, conification result holds for any
K\"ahler manifold endowed with a Hamiltonian
Killing vector field, see Theorem \ref{confThm}.  
This construction is new and 
may provide the needed conification procedure for the 
r-map. We will study this problem in the future. 

For the c-map we prove the existence of a canonical Killing vector field
satisfying the assumptions of Theorem \ref{hKconifThm}. In this way we 
can associate a 
family of (possibly indefinite) conical hyper-K\"ahler metrics and, hence, 
a family of quaternionic K\"ahler metrics to any
projective special K\"ahler manifold. In view of the results
of \cite{APP} Section 2.4, this  
family should contain the Ferrara-Sabharwal metric as well as 
the (locally defined) one-parameter deformation discussed in \cite{APSV}
and in the papers cited there. 
The parameter should be related to the choice of Hamiltonian 
function $f$, which is unique up to a constant.  
This will be the topic of future investigation. 

\noindent
{\bf Acknowledgments}
We thank Stefan Vandoren for discussions and 
for his notes about examples of the 
QK/HK correspondence. We also thank Malte Dyckmanns for 
useful comments. This work is part of a research 
project within the RTG 1670 ``Mathematics inspired by String Theory'',
funded by the Deutsche Forschungsgemeinschaft (DFG).  The work of T.M.\ is 
supported in part by STFC grant ST/J000493/1. T.M.\ thanks the 
Department of Mathematics and the SFB 676 for hospitality and support 
during several stages of this work.

\section{Conification of K\"ahler manifolds} 
\bd An {\cmssl almost pseudo-K\"ahler manifold} 
$(M,g,J)$ is a pseudo-Riemannian 
manifold $(M,g)$ endowed with a skew-symmetric almost complex structure
with closed fundamental form $\o := g(J\cdot ,\cdot )$. 
It is called a {\cmssl pseudo-K\"ahler manifold} if the almost
complex structure is integrable. In that case $\o$ is called the 
{\cmssl K\"ahler form}. 
\ed 

Let $(M, g , J, Z)$ be a pseudo-K\"ahler manifold 
endowed with a time-like or space-like Hamiltonian Killing vector field 
$Z$.  Let $-f$ be a corresponding Hamiltonian function, that is 
$df=-\o (Z,\cdot )$, where $\o$ is the
K\"ahler form. We will assume that
$f$ and $f_1:=f-\frac{1}{2}g(Z,Z)$ are nowhere vanishing. 
\bl \label{dhLemma} 
Let $Z$ be a Killing vector field on a pseudo-K\"ahler manifold
$(M,g,J)$ and put $h=\frac{g(Z,Z)}{2}$. Then 
\[ dh = \o (JD_ZZ ,\cdot ),\]
where $D$ is the Levi-Civita connection.  
In particular, $dh=-\o (Z,\cdot )$ holds if and only if 
$D_ZZ=JZ$. 
\el  
\pf 
\[ dh= g(DZ,Z)=-g(D_ZZ,\cdot ) =\o (JD_ZZ,\cdot ).\] 
\qed 

\noindent 
The Lemma implies that 
\be \label{f1Equ} df_1 =  d(f-h)= -\o (Z + JD_ZZ,\cdot )
= -g(J(Z + JD_ZZ),\cdot ).\ee

Let $\pi : P\ra M$ be an 
$S^1$-principal bundle
endowed with a principal connection $\eta$ with the
curvature $d\eta = \pi^*(\o -\frac{1}{2}d\b)$, where $\b=g(Z,\cdot )$. 
Notice that locally we can always assume $P=M\times S^1$
and $\eta = ds + \eta_M$, where $\eta_M\in \O^1(M)$ and 
$s$ is the angular coordinate on $S^1=\{ e^{is}|s\in \bR\}$. 
We will denote the fundamental vector field of $P$ by $X_P$.  
It coincides with the vertical coordinate vector field $\p_s$
in any local trivialisation of the principal bundle. 
We define a pseudo-Riemannian metric on $P$ by
\[ g_P := \frac{2}{f_1}\eta^2 + \pi^*g.\]

Next we consider $\hat{M} := P\times \bR$  with the coordinate
$t$ on the $\bR$-factor and the projection
$\hat{\pi} : \hat{M}\ra M$ defined as $\hat{\pi} (p,t)= \pi (p)$,
for all $(p,t)\in \hat{M}$. On $\hat{M}$ we introduce the 
following tensor fields. 
\bea
\xi  &:=&  \p_t\in \mathfrak{X}(\hat{M}),\\ 
\hat{g} &:=& e^{2t}(g_P + 2fdt^2+2\a dt)\in \G (S^2T^*\hat{M}),
\label{hatgEqu}\\
\theta &:=& e^{2t}(\eta + \frac{1}{2}\b )\in \O^1(\hat{M}),\\
\hat{\o} &:=& d\theta\in \O^2(\hat{M}), 
\eea 
where $\a := df$ and covariant tensor fields on $M$ and $P$ 
are identified with their
pullbacks to tensor fields on $\hat{M}$.  

\bd A {\cmssl conical} pseudo-Riemannian manifold $(M,g,\xi)$ is a 
pseudo-Riemannian manifold  $(M,g)$ endowed with a time-like or space-like
vector field $\xi$ such that $D\xi = \mathrm{Id}$. 
\ed

\bt \label{confThm} Given $(M, g , J, Z)$ as above, the tensor 
field $\hat{g}$ defines a 
pseudo-Riemannian metric such that 
$(\hat{M}, \hat{g}, \hat{J} := \hat{g}^{-1}\hat{\o}, \xi )$ is  a conical   
pseudo-K\"ahler manifold.   The induced 
CR-structure on the hypersurface $P\subset \hat{M}$ coincides 
with the horizontal distribution $T^hP$ for the connection $\eta$ and
$\pi : P \ra M$ is holomorphic, that is 
$d\pi \hat{J}=Jd\pi$ on $T^hP$. The  
projection $\hat{\pi} : \hat{M} \ra M$ is not holomorphic, since 
$\mathrm{ker}\, d\hat{\pi} = \mathrm{span}\{X_P,\xi\}$ 
is not $\hat{J}$-invariant. The metric $\hat{g}$
has signature $(2k+2,2\ell )$ if $f_1>0$ and 
 $(2k,2\ell +2)$ if $f_1<0$, where $(2k,2\ell )$
is the signature of the metric $g$. 
\et 

\pf 
It is clear that the restriction of $\hat{g}$ to the horizontal
distribution $T^hP=\mathrm{ker}\, \eta \subset TP$ is nondegenerate. 
Let us denote by $E$ the orthogonal complement of 
the $2$-dimensional distribution 
$\mathrm{span}\{\widetilde{Z},\widetilde{JZ}\}\subset T^hP$, where 
$\widetilde{X}\in \mathfrak{X}(P)$ stands for the horizontal lift  
of a vector field $X\in \mathfrak{X}(M)$. Since $g(Z,Z)\neq 0$, 
we see that  $E{\oplus} \bR\widetilde{Z}\subset T^hP$ is
nondegenerate. The orthogonal distribution in $\hat{M}$ is precisely  
\[\mathcal{D}=\mathrm{span}\{\widetilde{JZ},X_P,\xi\},\]   
as follows from $\a (Z)=df(Z)=-\o (Z,Z)=0$. 
The matrix representing the bilinear form
$\hat{g}|_\mathcal{D}$ with respect to the frame  
$(\widetilde{JZ},X_P,\xi)$ 
is given by 
\[ e^{2t}\left( \begin{array}{ccc} 
g(Z,Z)&0&-g(Z,Z)\\
0&\frac{2}{f_1}&0\\
-g(Z,Z)&0&2f
\end{array} \right)\, ,\]
which has the determinant $4e^{6t}g(Z,Z)\neq 0$.
This proves that $\hat{g}$ is nondegenerate. 
The signature of $\hat{g}$ can be easily read off 
from the above matrix. 
 
Let us prove next that the skew-symmetric 
endomorphism field $\hat{J} = \hat{g}^{-1}\hat{\o}$ is also nondegenerate. 
Calculating the differential of $\theta$ we obtain
\be \hat{\o} = 2e^{2t}dt\wedge (\eta + \frac{1}{2}\b ) + e^{2t}\o 
\label{hatoEqu}.\ee 
This formula immediately implies that $\hat{J}$ preserves
the distribution $E$ and $\hat{J}\wt{X}=\wt{JX}$ for all
$X\in \mathfrak{X}(M)$ which are perpendicular to $Z$ and $JZ$.  

\noindent 
{\bf Claim 1:} $\hat{J}$ preserves
the distribution   $T^hP$ and 
\be \hat{J}\wt{X}=\wt{JX}\quad \mbox{for all}\quad 
X\in \mathfrak{X}(M).\label{hatJEqu}\ee\\ 
It remains to  check \re{hatJEqu}, or equivalently, that 
$\hat{\omega}(\wt{X},\cdot ) = \hat{g}(\wt{JX},\cdot)$,  
for $X\in \{ Z, JZ\}$. Using the formulas \re{hatgEqu} 
and \re{hatoEqu}, we have 
\[ \hat{\omega}(\wt{Z},\cdot ) = -e^{2t}(\beta (Z)dt+\a ),\quad 
 \hat{\omega}(\wt{JZ},\cdot ) =-e^{2t}\beta,\]
\[ \hat{g}(\wt{JZ},\cdot ) = e^{2t}(-\a +\a (JZ)dt),\quad 
\hat{g}(\wt{Z},\cdot )=e^{2t}\beta .\] 
This proves Claim 1, since $\a (JZ)=-\beta (Z)$.

\noindent 
{\bf Claim 2:}
\[ \hat{J} X_P = -\frac{1}{f_1}(\widetilde{JZ} + \xi ).\] 
It suffices to check that 
\[ \hat{\o} (X_P,\cdot ) = -\frac{1}{f_1}\hat{g}(\widetilde{JZ}+\xi ,\cdot ).\]
Using \re{hatoEqu}, we see that the left-hand side is simply $-2e^{2t}dt$.
The right-hand side yields
\[ -\frac{e^{2t}}{f_1}(-\a  +
\a (JZ)dt + 2fdt + \a ) = -\frac{e^{2t}}{f_1}( -g(Z,Z) +2f)dt= -2e^{2t}dt.\]
This proves Claim 2.\\
{\bf Claim 3:} 
\be \label{decompEqu} T\hat{M} = T^hP \stackrel{\perp}{\oplus} \mathrm{span}
\{ X_P,\hat{J}X_P\}. \ee 
In view of Claim 2, is clear that $X_P\perp T^hP$ and 
$\hat{J}X_P\perp E\oplus \bR \wt{Z}$. Therefore
it suffices to show that $\hat{J}X_P$ is perpendicular
to $\wt{JZ}$. 
We compute 
\[ -f_1\hat{g}(\hat{J}X_P,\wt{JZ}) = \hat{g}(\wt{JZ} + \xi , \wt{JZ}) = 
e^{2t}(g(Z,Z) + \a (JZ)) = 0.
\] 
This proves Claim 3.\\ 
{\bf Claim 4:} The distributions $T^hP, \mathrm{span}
\{ X_P,\hat{J}X_P\} \subset T\hat{M}$ are nondegenerate and orthogonal with 
respect to $\hat{\o}$.\\
In fact, 
\be \hat{\o}|_{T^hP} = e^{2t}\o \label{restr_omegaEqu}\ee 
is nondegenerate and also 
\[ \hat{\o} (X_P,\hat{J}X_P)=  2e^{2t}\left( dt \wedge (\eta +\frac{1}{2}\b )
\right) (X_P,
\hat{J}X_P) = -2e^{2t}dt(\hat{J}X_P) =  \frac{2e^{2t}}{f_1} \neq 0.\] 
The $\hat{\o}$-orthogonality of the distributions follows from Claim 3
and the $\hat{J}$-invariance of $T^hP$. 
This proves Claim 4.\\
{\bf Claim 5:} $\hat{J}$ is an almost complex structure.\\ 
Recall that, by Claim 1, $\hat{J}|_{T^hP}$ 
corresponds to the complex structure 
$J$ by means of the identification $T^hP\cong 
\pi^*TM$. Therefore, it suffices to check that $\hat{J}$
squares to $-\mathrm{Id}$ on $\mathrm{span} \{ X_P,\hat{J}X_P\}$. 
Using Claim 1 and 2, we compute
\[ \hat{J}^2X_P= \frac{1}{f_1}(\wt{Z} - \hat{J}\xi ) .\]
So we need to check that 
\be X_P = -\frac{1}{f_1}(\wt{Z} -\hat{J}\xi ),
\ee
or, equivalently, 
that 
\[ \hat{g}(X_P, \cdot ) = -\frac{1}{f_1}(\hat{g} (\wt{Z},\cdot ) 
-\hat{\o}(\xi ,\cdot ) ). \]
The left-hand side is simply $\frac{2e^{2t}}{f_1}\eta$ and the 
right-hand side 
\[ -\frac{e^{2t}}{f_1}(g(Z,\cdot ) -2\eta -\b ) = \frac{2e^{2t}}{f_1}\eta 
.\]
This proves Claim 5.\\ 
So far we have proven that $\hat{J}$ is a skew-symmetric 
almost complex structure with closed fundamental form, in other words that 
$(\hat{M}, \hat{g}, \hat{J})$ is  an almost pseudo-K\"ahler manifold. 
Notice that Claim 1 
implies that the induced CR-structure on $P$ coincides with
the horizontal distribution and that $\pi : P \ra M$ is holomorphic. 
Claim 2 shows that $\hat{\pi} :\hat{M} \ra M$ is not holomorphic. 

Next we prove that $D\xi =\mathrm{Id}$. 
By the Koszul formula, we have
\begin{eqnarray*} 2\hat{g}(D_{X_1}\xi ,X_2) &=& 
 X_1\hat{g}(\xi ,X_2)+
\xi \hat{g}(X_1,X_2)-X_2\hat{g}(X_1,\xi )\\ 
&&+\hat{g}([X_1,\xi],X_2)-\hat{g}(X_1,[\xi ,X_2]) 
-\hat{g}(\xi ,[X_1,X_2])\end{eqnarray*}
for all vector fields $X_1, X_2$ on $\hat{M}$. 
If $X_1, X_2$ are horizontal lifts of commuting
vector fields on $M$, the right-hand side 
yields
\[ e^{2t}d\a (X_1,X_2) + 2\hat{g}(X_1,X_2) = 2\hat{g}(X_1,X_2).\]
Similarly, if $X_1, X_2\in \{X_P, \xi \}$, the right-hand side
is also 
\[ 2f(X_1e^{2t}dt(X_2) -  X_2e^{2t}dt(X_1))
+ 2\hat{g}(X_1,X_2) = 2\hat{g}(X_1,X_2).\]  
Next we consider the case where $X_1$ is a horizontal lift and
$X_2=\xi$. The Koszul formula gives again 
\[ 2\hat{g}(D_{X_1}\xi ,\xi ) = 2e^{2t}X_1f = 2\hat{g}(X_1,\xi ) \]
and, similarly, for $X_2=X_P$:   
\[  2\hat{g}(D_{X_1}\xi , X_P) = 0 = 2\hat{g}(X_1,X_P).\]
Next, let $X_2$ be a horizontal lift and $X_1=\xi$.
Then 
\[ 2\hat{g}(D_{\xi}\xi , X_2) = 2\xi e^{2t} \a (X_2)- 2e^{2t}X_2f
= 2 e^{2t}\a (X_2)= 2\hat{g}(\xi , X_2).\]  
Finally, for   $X_1=X_P$ we get
\[2\hat{g}(D_{X_P}\xi , X_2) = 0 = 2\hat{g}(X_P,X_2).\] 

Next we prove that $\hat{J}$ is integrable. 
In order to apply the Newlander-Nirenberg theorem, let us
first recall that the decomposition \re{decompEqu}
is $\hat{J}$-invariant, in virtue of Claim 4. 
Therefore, Claim 1 implies that 
\[ T^{1,0}_p\hat{M} = \mathrm{span} \{ (\tilde{X}-i\widetilde{JX})_p|
X\in \mathfrak{X}(M)\} \oplus \bC (X_P -i \hat{J}X_P)_p,\]
for all $p\in \hat{M}$.  By the integrability
of the complex structure on $M$, we know that
for all $X, Y\in  \mathfrak{X}(M)$ there exists 
$W\in  \mathfrak{X}(M)$ such that 
\[ [X-iJX,Y-iJY] = W -iJW.\] 
Therefore, 
\begin{eqnarray*} [ \tilde{X}-i\widetilde{JX},  \tilde{Y}-i\widetilde{JY}] &=&
 \tilde{W}-i\widetilde{JW} - d\eta (\tilde{X}-i\widetilde{JX},  
\tilde{Y}-i\widetilde{JY})X_P\\
&=&
 \tilde{W}-i\widetilde{JW} - (\o -\frac{1}{2}d\beta ) (X-iJX,Y-iJY)X_P.
\end{eqnarray*}
Here we have used, the well known fact that the vertical part 
of the commutator of two horizontal vector fields on a principal 
bundle with connection is given by minus the curvature. 
We claim that not only $\omega$ but also  $d\beta$ 
is of type $(1,1)$, which finally implies 
$[ \tilde{X}-i\widetilde{JX},  \tilde{Y}-i\widetilde{JY}] = 
\tilde{W}-i\widetilde{JW}$. In fact, 
\[ d\beta = d g(Z,\cdot ) = -d\o (JZ ,\cdot ) = 
-\mathcal{L}_{JZ}\o\] 
is the Lie derivative of a form of type $(1,1)$ with respect to a
holomorphic (thus type-preserving) vector field. 
Finally, with the help of Claim 2, 
for $X\in  \mathfrak{X}(M)$, we compute   
\begin{eqnarray*} && [ \tilde{X}-i\widetilde{JX},X_P -i \hat{J}X_P] = 
-i[ \tilde{X}-i\widetilde{JX},\hat{J}X_P]\\
&=&
id(\frac{1}{f_1})(X-iJX)(\widetilde{JZ}+\xi ) 
+\frac{i}{f_1}[ \tilde{X}-i\widetilde{JX},
\widetilde{JZ}]\\
&=&i\frac{df_1}{f_1}(X-iJX)\hat{J}X_P-\frac{i}{f_1}d\eta (X-iJX,JZ)X_P +
\frac{i}{f_1}[ X-iJX,
JZ]^{\wt{}}.
\end{eqnarray*}
Notice that the last term is the horizontal lift of a vector
field of type $(1,0)$. In fact, it suffices to observe that the 
Lie derivative with respect to the holomorphic vector field 
$JZ$ preserves the type. The remaining part is of
type $(1,0)$ if and only if 
\be d\eta (X-iJX,JZ) = -idf_1(X-iJX). \label{auxiliaryEqu} \ee 
for all $X\in TM$. 
Now 
\begin{eqnarray*} d\eta (\cdot , JZ) &=& 
g(Z,\cdot ) + \frac{1}{2}d\b (JZ,\cdot ),\\ 
d\b (JZ,\cdot ) &=& \mathcal{L}_{JZ}\b = g(D_Z(JZ),\cdot ) +
g(Z,D(JZ))\\
g(Z,D(JZ)) &=& - g(JZ,DZ) = g(D_{JZ}Z, \cdot ) = g(D_Z(JZ),\cdot )
= g(JD_ZZ,\cdot )
\end{eqnarray*}
Therefore,  
\[ d\eta (\cdot , JZ) = g(Z+JD_ZZ,\cdot ).\]
Comparing with \re{f1Equ} we see that 
that \re{auxiliaryEqu} is equivalent to 
\begin{eqnarray*} &&g(Z + JD_ZZ,X-iJX) = ig(J(Z+JD_ZZ),X-iJX)\\ 
&=& g(J(Z+JD_ZZ),J(X-iJX))= g(Z + JD_ZZ,X-iJX) ,
\end{eqnarray*}
which is always satisfied. 
\qed \\

\bd Let $(M,g)$ be any pseudo-Riemannian manifold. Then 
$C_\pm (M):=(\mathbb{R}^{>0}\times M, \pm dr^2+r^2g)$ is called the 
{\cmssl space-like or time-like cone} over $(M,g)$, respectively.  
\ed 
The vector field $\xi = r\p_r$ defines on $C_\pm (M)$ the structure
of a conical pseudo-Riemannian manifold and any conical  pseudo-Riemannian 
manifold is locally isomorphic to a  space-like or time-like cone.  

\bd A {\cmssl pseudo-Sasakian structure} on a pseudo-Riemannian manifold
$(M,g)$ is a unit Killing vector field $Z$ such that 
$J:= DZ|_{Z^\perp}$ defines an integrable CR-structure $H=Z^\perp \subset TM$ 
with the Levi form $2g$. The {\cmssl Levi-form} is the symmetric 
bilinear form $L$ on $H$ defined by $L(X,Y)= 
\frac{g(Z,[JX,Y])}{g(Z,Z)}$. 
\ed   
It is well known that $(M,g)$ admits a 
space-like or time-like pseudo-Sasakian structure $Z$ if and only 
if the space-like  or time-like cone over $(M,g)$
admits a K\"ahler structure $\hat{J}$ compatible
with the cone metric.  

\begin{Ex} In Theorem \ref{confThm} we have assumed that $Z$ is 
nowhere light-like. However, one can verify that 
the construction remains meaningful if we put $Z=0$. 
Taking $Z=0$ and $f=const=c\neq 0$ in the construction
of Theorem \ref{confThm}, yields a conical pseudo-K\"ahler manifold  
$(\hat{M}, \hat{g}, \hat{J}= \hat{g}^{-1}\hat{\o}, \xi )$. It is 
precisely the space-like ($c>0$) or time-like ($c<0$) cone over 
$(P,\frac{1}{2|c|}g_P)$, where $r=\sqrt{2|c|}e^t$. The 
unit Killing vector field $\zeta := |c| X_P$ defines a 
pseudo-Sasakian structure on $(P,\frac{1}{2|c|}g_P)$.
Notice that   $(P,\frac{1}{2|c|}g_P)$ is a pseudo-Riemannian
submersion over the pseudo-K\"ahler manifold $(M,\frac{1}{2|c|}g)$.   
In particular, we can take $f=\pm \frac{1}{2}$ and $r=e^t$,
which yields $(\hat{M}, \hat{g})$ as the space-like or time-like cone over 
the pseudo-Sasaki manifold $(P,g_P,\zeta=\frac{1}{2}X_P)$ and the latter 
fibers as a pseudo-Riemannian submersion over $(M,g)$. Alternatively,
we may take $c=\pm 1$, for which $X_P$ is the Sasaki structure. In that case
$(\hat{M}, \hat{g})$ is the space-like or time-like cone over 
the pseudo-Sasaki manifold $(P,\frac{1}{2}g_P= \pm \eta^2 + \frac{1}{2}g,X_P)$ 
and the latter fibers as a pseudo-Riemannian submersion over $(M,\frac{1}{2}g)$.
\end{Ex}

\section{Conification of hyper-K\"ahler manifolds}\label{HKSec}
Let  $(M, g , J_1, J_2, J_3)$ be a pseudo-hyper-K\"ahler manifold 
with the three K\"ahler forms $\o_\a := gJ_\a:= g(J_\a \cdot ,\cdot )$,
$\a =1,2,3$. We will assume that $Z$ is a time-like or space-like 
Killing vector field and that $f$ is a nowhere vanishing function 
such that $df = -\o_1Z:=-\o_1(Z,\cdot )$. Following the notation of the 
previous section, we put $f_1 := f- h$, where $h :=\frac{g(Z,Z)}{2}$. 
We will also assume that $f_1$ is nowhere zero.  
Applying Theorem \ref{confThm} to the pseudo-K\"ahler manifold 
$(M, g , J_1)$ endowed with the $\o_1$-Hamiltonian Killing vector field $Z$, 
we obtain the principal bundle $\pi : P\ra M$ with 
the connection $\eta$ and the pseudo-Riemannian metric $g_P$
such that $\hat{M}_1:=P\times \bR$ is endowed 
with the structure of a conical pseudo-K\"ahler manifold.  
Our aim is to construct a conical pseudo-hyper-K\"ahler manifold
$(\hat{M},\hat{g},\hat{J}_1,\hat{J}_2,\hat{J}_3,\xi )$ such that 
$\hat{M}_1\subset \hat{M}$ with the conical pseudo-K\"ahler structure induced
by $(\hat{g},\hat{J}_1,\xi )$.  
As a first step, we define the vector field 
\[ Z_1:= \tilde{Z} +f_1X_P\] 
and the one-forms  
\begin{eqnarray}\label{thetaEqu}
\theta_1^P &:=& \eta +\frac{1}{2}gZ\nonumber \\ 
\theta_2^P  &:=& \frac{1}{2}\o_3Z\nonumber \\ 
\theta_3^P  &:=& -\frac{1}{2}\o_2Z 
\end{eqnarray}
on $P$.
We consider $\theta_\a := f^{-1}\theta_\a^P$ as the 
components of a one-form $\theta := \sum_\a \theta_\a i_\a$ with values in the 
imaginary quaternions, where $(i_1,i_2,i_3)=(i,j,k)$. 
Then we extend  $\theta$ to a one-form $\tilde{\theta}$ on  
$\tilde{M} := \bH^*\times P \supset \{ 1\}\times P\cong P$ by
\[ \tilde{\theta}_\a (q,p) := \varphi_\a (q) + (\mathrm{Ad}_q\theta (p))_\a ,\quad 
(q,p)\in \tilde{M},\]
where $\varphi = \varphi_0 + \sum_\a \varphi_\a i_\a$ is 
the right-invariant Maurer-Cartan form of 
$\bH^*$ and $\mathrm{Ad}_qx = qxq^{-1} = x_0 + \sum_\a (\mathrm{Ad}_qx)_\a i_\a$  for all 
$x= x_0 + \sum x_\a i_\a \in \bH$. Notice that
\[ \varphi_a(e_b) = \d_{ab},\]
where $(e_0,\ldots ,e_3)$, is the right-invariant frame  of $\bH^*$ 
which coincides with the standard basis of 
$\bH=\mathrm{Lie}(\bH^*)$ at $q=1$. 
Next we define 
\[ \tilde{\o}_\a := d (\rho^2\tilde{\theta}_\a^P),\] 
where $\tilde{\theta}_\a^P:= f\tilde{\theta}_\a$ and $\rho := |q|$. 
Let us denote by $e_1^L$ the left-invariant vector field on
$\bH^*$ which coincides with $e_1$ at $q=1$ and by $\hat{M}$ 
the space of integral curves of the vector field $V_1:=e_1^L-Z_1$. 
We will assume that the quotient map $\tilde{\pi} : \tilde{M}\ra \hat{M}$
is a submersion onto a Hausdorff manifold. (Locally this is
always the case, since the vector field has no zeroes.)
\bt \label{hKconifThm} 
Let  $(M, g , J_1, J_2, J_3, Z)$ be a pseudo-hyper-K\"ahler manifold
endowed with a Killing vector field $Z$ satisfying the above assumptions 
and $\mathcal{L}_ZJ_2= -2J_3$.  
Then there exists a pseudo-hyper-K\"ahler structure 
$(\hat{g},\hat{J}_1 ,\hat{J}_2 ,\hat{J}_3)$
on $\hat{M}$ with exact K\"ahler forms $\hat{\o}_{\a}$ determined by
\be \label{hatoaEqu}\tilde{\pi}^*\hat{\o}_{\a}=\tilde{\o}_\a .\ee  
The vector field $\r \partial_\r$ on $\tilde{M}$ projects to a 
vector field $\xi$ on $\hat{M}$ such that  $(\hat{M},
\hat{g},\hat{J}_1 ,\hat{J}_2 ,\hat{J}_3,\xi)$ is a 
conical pseudo-hyper-K\"ahler manifold. The signature
of the metric $\hat{g}$ is $(4k,4\ell +4)$
if  $f_1<0$ and $(4k+4,4\ell )$ if $f_1>0$, where
$(4k,4\ell )$ is the signature of the metric $g$. 
\et 

\pf We first show that the one-forms $\tilde{\theta}_\a$ on $\tilde{M}$
induce one-forms $\hat{\theta}_\a$ on the quotient $\hat{M}$. 
\bl There exist one-forms
$\hat{\theta}_\a$ on $\hat{M}$ such that $\tilde{\theta}_\a= 
\tilde{\pi}^*\hat{\theta}_\a$.
\el

\pf 
Let us first observe that 
the above definitions imply that $\theta (Z_1)=i_1=i$. 
To compute $\varphi (e_1^L)$, we use the equivariance
of the right-invariant Maurer-Cartan form with respect to left-multiplication: 
\[ \varphi (dL_qv) = \mathrm{Ad}_q\varphi (v),\]
for all $q\in \bH^*$, $v\in T\bH^*$.  Using that $\varphi (v)=v$ for all 
$v\in T_e\bH^*$, we conclude that  
\[ \varphi (e_\a^L)= \varphi (dL_qi_\a) = \mathrm{Ad}_q(i_\a ).\] 
Combining these facts, we get 
\[ \tilde{\theta}(V_1) = \varphi (e_1^L) - \mathrm{Ad}_q(\theta (Z_1))
= \mathrm{Ad}_qi-\mathrm{Ad}_qi=0.\]
This shows that the one-forms  $\tilde{\theta}_\a$ on $\tilde{M}$
are horizontal with respect to the projection $\tilde{\pi} : 
\tilde{M}\ra \hat{M}$.  To prove the lemma, it now suffices to   
check that $\mathcal{L}_{V_1}\tilde{\theta}=0$. 
First of all, the right-invariance
of $\varphi$ implies the invariance under any left-invariant vector 
field. So  $\mathcal{L}_{V_1}\varphi = \mathcal{L}_{e_1^L}\varphi=0$
and we are left with 
\be \label{V1Equ} \mathcal{L}_{V_1}\tilde{\theta}=
\mathcal{L}_{e_1^L}\mathrm{Ad}_q\theta  - 
\mathrm{Ad}_q\mathcal{L}_{Z_1}\theta.\ee 
The first term is easily computed as follows:
\[ \mathcal{L}_{e_1^L}\mathrm{Ad}_qx =
\frac{d}{dt}\left|_{t=0}q\exp(ti)x\exp(-ti)q^{-1}\right.=  
\mathrm{Ad}_q[i,x],\]
for all $x\in \bH$.  This shows that 
\be \label{e1LEqu} 
\mathcal{L}_{e_1^L}\mathrm{Ad}_q\theta = \mathrm{Ad}_q[i,\theta].\ee
For the computation of $\mathcal{L}_{Z_1}\theta$ we first remark
that $\mathcal{L}_{X_P}\theta=0$, such that  $\mathcal{L}_{f_1X_P}\theta =
f^{-1}df_1\theta(X_P)= i f^{-1}df_1$ and $\mathcal{L}_{Z_1}\theta
= \mathcal{L}_{\tilde{Z}}\theta +if^{-1}df_1$. We compute each component 
$\mathcal{L}_{\tilde{Z}}\theta_\a$. Using that $Z$ is Killing, we get  
\[ \mathcal{L}_{\tilde{Z}}\theta_1^P = \mathcal{L}_{\tilde{Z}}\eta 
+ \frac{1}{2}\mathcal{L}_{Z}(gZ)= \mathcal{L}_{\tilde{Z}}\eta
= \o_1Z -\frac{1}{2}d(gZ)Z=  \o_1Z+dh=-df_1.
\]
Here we used that $d(gZ)Z=\mathcal{L}_{Z}(gZ)-d(g(Z,Z))=0-2dh=-2dh$. 
The hypothesis $\mathcal{L}_ZJ_2= -2J_3$ on the $\o_1$-Hamiltonian 
Killing vector field
$Z$ immediately implies 
\[ \mathcal{L}_{\tilde{Z}}\theta_2^P= \frac{1}{2}\mathcal{L}_{Z}\o_3Z =
\frac{1}{2}\mathcal{L}_{Z}\o_1J_2Z = - \o_1J_3Z=\o_2Z=-2\theta_3^P \]
and, similarly, $\mathcal{L}_{\tilde{Z}}\theta_3^P=2\theta_2^P$. 
(Notice that $\mathcal{L}_ZJ_2= -2J_3$ implies $\mathcal{L}_ZJ_3= 2J_2$,
because $Z$ is $J_1$-holomorphic.)   
Since $\mathcal{L}_{Z_1}f=0$, this shows that $\mathcal{L}_{Z_1}\theta_1=
\mathcal{L}_{Z_1}(f^{-1}\theta_1^P)=0$, $\mathcal{L}_{Z_1}\theta_2=-2\theta_3$
and $\mathcal{L}_{Z_1}\theta_3=2\theta_2$. Summarising, we have 
\be  \label{Z1Equ} \mathcal{L}_{Z_1}\theta = [i,\theta ].\ee 
The equations \re{V1Equ}, \re{e1LEqu} and \re{Z1Equ} show that 
\[  \mathcal{L}_{V_1}\tilde{\theta}= \mathrm{Ad}_q[i,\theta ]-
\mathrm{Ad}_q[i,\theta ]=0.\] 
\qed 

Since $\mathcal{L}_{V_1}\rho=\mathcal{L}_{V_1}f=0$, the functions 
$f$ and $\rho$ are well defined on the quotient $\hat{M}$. 
Therefore, the lemma shows that 
\[ \hat{\o}_\a := d(\rho^2f\hat{\theta}_\a )\]
are two-forms on $\hat{M}$, which satisfy \re{hatoaEqu}. 
To prove that  the triplet $(\hat{\o}_\a)$ defines a pseudo-hyper-K\"ahler
structure we will prove that the $\hat{\o}_\a$ are nondegenerate 
such that we can consider the nondegenerate endomorphisms  $\hat{J}_\a$ defined
by 
\be \label{hatJaEqu} \hat{\o}_\a  \hat{J}_\b = \hat{\o}_\g\ee 
for any cyclic permutation $(\a , \b , \g )$
of $(1,2,3)$. In the following $(\a , \b , \g )$ will be always a 
cyclic permutation. We have to show that $(\hat{J}_\a )$ is an almost 
hyper-complex
structure. The integrability then 
follows from the closure of the $\hat{\o}_\a$, 
in virtue of the Hitchin lemma. The pseudo-hyper-K\"ahler metric 
is then given by $\hat{g}=-\hat{\o}_\a \hat{J}_\a$. Notice that this 
expression is independent of $\a$, due to the relations $\hat{J}_\a \hat{J}_\b = 
\hat{J}_\g$. The skew-symmetry of $\hat{J}_\b$ with respect to $\hat{g}$ 
follows from the symmetry of $\hat{J}_\b$ with 
respect to $\hat{\o}_\a$ (a  
consequence of  \re{hatJaEqu}) and the relation $\hat{J}_\a \hat{J}_\b = 
-\hat{J}_\b \hat{J}_\a$.  The symmetry 
of $\hat{g}$ is then obtained from
$\hat{J}_\a^2=-\mathrm{Id}$. The nondegeneracy of $\hat{g}$ is a consequence 
of that of $\hat{\o}_\a$ and $\hat{J}_\a$. 
\bl \label{tildooLemma} The two-forms $\tilde{\o}_\a$ on $\tilde{M}$ are given by:
\be \label{tildooEqu}\tilde{\o}_\a = 2f\r^2(\varphi_0\wedge \varphi_\a +\varphi_\b\wedge 
\varphi_\g  + \varphi_0\wedge \theta_\a'-\theta_0\wedge
\varphi_\a + \varphi_\b\wedge \theta_\g'-\varphi_\g\wedge \theta_\b')
+\r^2\o' ,\ee
where $\theta_0:=-\frac{1}{2}f^{-1}df$, 
$\theta' := \mathrm{Ad}_q\theta$, 
$\o := \sum  \o_\a i_\a$  
and $\o' =\mathrm{Ad}_q\o$, 
\el 

\pf We first calculate the differential of $\tilde{\theta}^P = f\tilde{\theta} 
= f\varphi + \mathrm{Ad}_q\theta^P$, where $\theta^P=f\theta$. 
Using the Maurer Cartan equation 
\[ d\varphi_\a = 2\varphi_\b\wedge 
\varphi_\g,\] 
we obtain  
\[ d(f\varphi_\a )= 2f(-\theta_0\wedge \varphi_\a + \varphi_\b\wedge 
\varphi_\g ).\] 
Using the fact that $\varphi = dqq^{-1}$, 
we see that
\begin{eqnarray*} d(\mathrm{Ad}_q\theta^P) &=& 
dq\wedge \theta^Pq^{-1}+
qd\theta^Pq^{-1} - q\theta^P\wedge d(q^{-1})\\
&=&\varphi\wedge \mathrm{Ad}_q\theta^P 
+\mathrm{Ad}_q(\theta^P) \wedge \varphi + \mathrm{Ad}_qd\theta^P\\
&=& f(\varphi\wedge \theta'
+\theta' \wedge \varphi) + \mathrm{Ad}_qd\theta^P .
\end{eqnarray*} 
The components are given by
\begin{eqnarray*}  d(\mathrm{Ad}_q\theta^P)_\a &=& 
f(\varphi_0\wedge \theta_\a' + 
\varphi_\b\wedge \theta_\g' - \varphi_\g\wedge \theta_\b' + 
\theta_\b' \wedge \varphi_\g -\theta_\g' \wedge \varphi_\b + \theta_\a'\wedge 
\varphi_0 ) + 
(\mathrm{Ad}_qd\theta^P)_\a \\
&=& 2f(\varphi_\b\wedge \theta_\g' - \varphi_\g\wedge \theta_\b' )  + 
(\mathrm{Ad}_qd\theta^P)_\a .
\end{eqnarray*}  
Using $\varphi_0=\r^{-1}d\r$ and the above equations, we get 
\begin{eqnarray*} \tilde{\o}_\a &=& d(\r^2\tilde{\theta}^P_\a ) 
= 2 \r^2\varphi_0\wedge 
\tilde{\theta}^P_\a + \r^2d\tilde{\theta}^P_\a= 2 f\r^2\varphi_0\wedge 
(\varphi_\a + \theta_\a') + \r^2d\tilde{\theta}^P_\a\\
&=& 2 f\r^2(\varphi_0\wedge \varphi_\a +\varphi_\b\wedge 
\varphi_\g  + \varphi_0\wedge \theta_\a'-\theta_0\wedge
\varphi_\a + \varphi_\b\wedge \theta_\g'-\varphi_\g\wedge \theta_\b') + 
\r^2(\mathrm{Ad}_qd\theta^P)_\a  .
\end{eqnarray*}
Finally, we claim that 
\[ d\theta^P = \o ,\]
which implies the lemma. 
In fact, 
\begin{eqnarray*} d\theta_1^P &=& d\eta + \frac{1}{2}d(gZ)= \o_1\\
d\theta_2^P &=& \frac{1}{2}d(\o_3Z)= \frac{1}{2}\mathcal{L}_Z\o_3 =
\o_2 \\
d\theta_3^P &=& -\frac{1}{2}d(\o_2Z)= -\frac{1}{2}\mathcal{L}_Z\o_2 =
\o_3. \\
\end{eqnarray*}
\qed 

Next we will show that the two-forms $ \tilde{\o}_\a$ computed
in Lemma \ref{tildooLemma} are nondegenerate on any 
distribution complementary to $\bR V_1\subset T\tilde{M}$.
Let us denote by $\mathcal{D}_1, \mathcal{D}_2\subset T\tilde{M}$ 
the distributions which correspond to the factors of the product 
$\tilde{M} =\bH^*\times P$. The second distribution can be 
decomposed as 
\[   \mathcal{D}_2 = \bR X_P \stackrel{\perp}{\oplus}  \mathcal{D}_2^h,\quad
\mathcal{D}_2^h =  E \stackrel{\perp}{\oplus} E',\quad E' := 
\mathrm{span}\{ \tilde{Z}, \wt{J_1Z}, \wt{J_2Z}, 
\wt{J_3Z}\}
,\] 
with respect to the metric $g_P$ on the leaves $\{ q\}\times P\cong P$
of  $\mathcal{D}_2$. Notice that $\mathcal{D}_2^h|_{(q,p)}\cong T_p^hP$. 
We will study the restriction of $\tilde{\o}_\a$ to 
the distribution $\mathcal{D}_1\oplus 
\mathcal{D}_2^h$, which is complementary to $\bR V_1$. 
{}From  \re{tildooEqu} we first see that the distributions $E$ and
$\mathcal{D} := \mathcal{D}_1\oplus E'$ are orthogonal with 
respect to $\tilde{\o}_\a$. Furthermore, 
\[ \tilde{\o}_\a|_E = \r^2\o_\a'|_E = \r^2 gJ_\a'|_E,\quad J_\a' 
= \sum A_{\a \b}J_\b,
\]
where $(A_{\a \b})\in \SO(3)$ is the matrix $A=A(q)$ representing
$\mathrm{Ad}_q|_{\mathrm{Im}\,\bH }$ in the basis $(i_1,i_2,i_3)$. This  
shows that  $\tilde{\o}_\a|_E$ is nondegenerate and that 
\[ \tilde{\o}_\a J_\b'|_E = \tilde{\o}_\g|_E.\] 
Notice that $(J_\a')$ coincides with the hyper-complex structure
$(J_\a)$ of $M$ up to a rotation, which depends on $q$.  
It remains to analyse  $\tilde{\o}_\a$ on the 
eight-dimensional distribution $\mathcal{D} = \mathcal{D}_1\oplus E'$.
We put 
\[ W_0 := J_1Z,\quad W_\a := -J_\a W_0\]
and 
\[ W_\a'=\sum_{\b =1}^3A_{\a \b}W_\b.\] 
{}From \re{thetaEqu} one can check that 
\[ \theta_a^P (\tilde{W}_b ) = h\d_{ab},\quad a,b\in \{0,\ldots ,3\},\]
where $\theta_0^P:=f\theta_0=-\frac{1}{2}df$. 
This implies that 
\[ \o_\a|_{E'} = -2h^{-1}(\theta_0^P\wedge \theta_\a^P -\theta_\b^P 
\wedge \theta_\g^P)= 
-2f^2h^{-1}(\theta_0\wedge \theta_\a -\theta_\b \wedge \theta_\g )\]
and, hence, 
\[ \tilde{\o}_\a|_{E'}=-2f^2h^{-1}\r^2(\theta_0\wedge \theta_\a' -
\theta_\b' \wedge \theta_\g' ).\] 
Now we can write the matrix $\mathcal{M}(\tilde{\o}_\a )$
which represents  $\tilde{\o}_\a|_{\mathcal{D}}$
in the basis 
$(e_0,e_\a,e_\b,e_\g,$ $\tilde{W}_0,\tilde{W}_\a' ,\tilde{W}_\b' ,
\tilde{W}_\g' )$: 
\be \label{tildoD} \mathcal{M}(\tilde{\o}_\a ) = 2\r^2\left( \begin{array}{rr} f\left( \begin{array}{rr} J&0\\
0&J
    \end{array}\right) 
    &h\left(\begin{array}{rr} I&0\\
0&J
    \end{array}\right) \\
-h\left( \begin{array}{rr} I&0\\
0&-J
    \end{array}\right) &-h\left( \begin{array}{rr} J&0\\
0&-J
    \end{array}\right) 
  \end{array}
\right) ,\ee
where
\[ I= \left( \begin{array}{rr} 0&1\\
1&0
    \end{array}\right) ,\quad J= \left( \begin{array}{rr} 0&1\\
-1&0
    \end{array}\right) .\]  

The invertibility of this matrix follows from the 
assumption $f_1=f-h\neq 0$. This proves that
the two-forms $\tilde{\o}_\a$ are nondegenerate
on any complement of $\bR V_1$, which implies the 
nondegeneracy of the induced forms $\hat{\o}_\a$ on $\hat{M}$.  
Now we compute the three endomorphisms 
\[ \tilde{J}_\a = \tilde{J}_\a|_{\mathcal{D}}\oplus \tilde{J}_\a|_E = 
\tilde{J}_\a|_{\mathcal{D}}\oplus J_\a'|_E\] 
of $\mathcal{D}\oplus E\cong T\tilde{M}/\bR V_1$ defined by
\[ \tilde{\o}_\a \tilde{J}_\b = \tilde{\o}_\g|_{\mathcal{D}\oplus E}.\] 
Under the projection $\tilde{M} \ra \hat{M}$, 
they correspond to the three endomorphism fields $\hat{J}_\a$ on $\hat{M}$
such that 
\[ \hat{\o}_\a \hat{J}_\b = \hat{\o}_\g .\] 
Using the expression \re{tildoD} one can check that the matrix 
representing $\tilde{J}_\a|_{\mathcal{D}}$ in the basis 
$(e_0,e_\a,e_\b,e_\g,\tilde{W}_0,\tilde{W}_\a' ,\tilde{W}_\b' ,
\tilde{W}_\g' )$ is given by 
\be \label{MJalpha} \mathcal{M}(\tilde{J}_\a ) = 
\left( \begin{array}{rrrr} -J&0&0&0\\
0&-J&0&0\\
0&0&J&0\\
0&0&0&-J
\end{array}
\right) ,
 \ee 
or, equivalently, 
\[ \begin{array}{ccc}
\tilde{J}_\a e_0 &=& e_\a\\ 
\tilde{J}_\a e_\a &=& -e_0\\
\tilde{J}_\a e_\b &=& e_\g\\
\tilde{J}_\a e_\g  &=& -e_\b
\end{array}
\quad\quad\quad 
\begin{array}{ccc}
\tilde{J}_\a\tilde{W}_0 &=& -\tilde{W}_\a'\\
\tilde{J}_\a\tilde{W}_\a' &=& \tilde{W}_0\\
\tilde{J}_\a\tilde{W}_\b' &=& \tilde{W}_\g'\\
\tilde{J}_\a\tilde{W}_\g' &=& -\tilde{W}_\b'.\\
\end{array}
\]
This proves that the $\tilde{J}_\a$ satisfy the standard 
quaternionic relations and that $\tilde{J}_\a|_{\mathcal{D}_2^h}$
corresponds to $J_\a'$ under the isomorphism $\mathcal{D}_2^h|_{(q,p)}\cong 
T^h_pP\cong T_{\pi (p)}M$, $(q,p)\in \tilde{M}$. Therefore, we have proven
that the three  symplectic forms  $\hat{\o}_\a$ on $\hat{M}$ define
a hyper-K\"ahler structure $(\hat{M},\hat{g},\hat{J}_\a, \a= 1,2,3)$.   

Next we calculate the explicit expression for the pseudo-hyper-K\"ahler 
metric $\hat{g}$. It amounts to calculating the metric 
\[ \tilde{g} := -\tilde{\o}_\a \tilde{J}_\a ,\]
which is defined on the codimension one distribution 
$\mathcal{D}\oplus E \subset T\tilde{M}$.  
\bp 
\[ \tilde{g} = \tilde{g}|_{\mathcal{D}} \oplus \tilde{g}_E,\quad 
\tilde{g}|_{\mathcal{D}} = 2\r^2 f (\sum_{a=0}^3\varphi_a^2 
+  h^{-1}\sum_{a=0}^3 (\theta_a' )^2 - 2\varphi_0\theta_0
+2\sum_{\a =1}^3 \varphi_\a \theta_\a') ,
\quad \tilde{g}_E= \r^2 g|_E.\]
\ep 

\pf
It suffices to calculate the matrix $\mathcal{M}(\tilde{g})$ which represents 
$\tilde{g} =  
-\tilde{\o}_\a \tilde{J}_\a $ 
in the basis $(e_0,e_\a,e_\b,e_\g,\tilde{W}_0,\tilde{W}_\a' ,\tilde{W}_\b' ,
\tilde{W}_\g' )$ of $\mathcal{D}$. In view of \re{tildoD} and \re{MJalpha},  it 
is given by
\be \label{tildegmatrixEqu} \mathcal{M}(\tilde{g}) =  
-\mathcal{M}(\tilde{J}_\a )^t\mathcal{M}(\tilde{\o}_\a ) = 
\mathcal{M}(\tilde{\o}_\a )  \mathcal{M}(\tilde{J}_\a )  = 
2\r^2 \left( \begin{array}{llllllll}f&&&&-h&&&\\
&f&&&&h&&\\
&&f&&&&h&\\
&&&f&&&&h\\
-h&&&&h&&&\\
&h&&&&h&&\\
&&h&&&&h&\\
&&&h&&&&h
\end{array}\right),
\ee 
where only the nonzero entries are written. This proves the above formula 
for $\tilde{g}_\mathcal{D}$, 
since $\theta_a'(\tilde{W}_b')=f^{-1}h\delta_{ab}$.  
\qed 

\noindent 
Let us now extend $\tilde{g}$ from a metric defined on the distribution 
$\mathcal{D}\oplus E =\mathcal{D}_1 \oplus\mathcal{D}_2^h 
\subset T\tilde{M}$ to a 
pseudo-Riemannian metric on
$\tilde{M}$ such that $V_1$ is perpendicular to $\mathcal{D}\oplus E$. 
Then $\tilde{\pi} : (\tilde{M},\tilde{g}) \ra
(\hat{M},\hat{g})$ is a pseudo-Riemannian submersion and we can
calculate the covariant derivative of
$\xi = \tilde{\pi}_\ast e_0$ by calculating $\tilde{g}(D_Xe_0,Y)$ 
for all $X, Y\in \mathcal{D}\oplus E$. In order to show that
$D\xi =\mathrm{Id}$, we have to check that
$\tilde{g}(D_Xe_0,Y)= \tilde{g}(X,Y)$. 
Using the Koszul formula and the commutator relations of the vector fields 
$e_a$, we obtain 
\[ 2\tilde{g}(D_{e_a}e_0,e_b) =  2f(\d_{0b}e_a + \d_{ab}e_0  -\d_{0a}e_b)\r^2 =
2f\d_{ab}e_0\r^2=4f\r^2\d_{ab}=2\tilde{g}(e_a,e_b),\]
for all $a,b\in \{0,\ldots ,3\}$.
Let us next observe that
\be \label{ge0Equ} 
\tilde{g}(e_0,\cdot )|_{\mathcal{D}_2}= -2\r^2 \theta_0^p,\ee  
as follows from $\tilde{g}(e_0,\cdot )|_{\mathcal{D}_2^h}= -2\r^2 \theta_0^p$ 
and $\tilde{g}(e_0,X_P) = \tilde{g}(e_0,\frac{e_1^L-\tilde{Z}}{f_1})= 
0$. 
Now let $X,Y\in \G (T^hP) \subset 
\G (\mathcal{D}_2^h)$ be horizontal lifts of vector fields in 
$M$.  Then using \re{ge0Equ},  $[e_0,X]= [e_0,Y]=0$ and   
$d\theta_0^P=0$  we get 
\begin{eqnarray*} 2\tilde{g}(D_Xe_0,Y) &=& X\tilde{g}(e_0,Y) + e_0\tilde{g}(X,Y)
- Y\tilde{g}(e_0,X) - \tilde{g}(e_0,[X,Y])\\
 &=& 2\tilde{g}(X,Y)
-2\r^2 (X \theta_0^P(Y) - Y\theta_0^P(X) -\theta_0^P([X,Y]))\\
&=&  2\tilde{g}(X,Y) -2\r^2 d\theta_0^P (X,Y) = 2\tilde{g}(X,Y)  .
\end{eqnarray*}

To compute $\tilde{g}(D_{e_a}e_0,X)$, we observe that $[e_0,e_a]=
[e_0,X]=[e_a,X]=0$, such that
\begin{eqnarray*}  2\tilde{g}(D_{e_a}e_0,X) &=& 
e_a\tilde{g}(e_0,X)+ e_0\tilde{g}(e_a,X)
- X\tilde{g}(e_0,e_a)\\ 
&=& 2\tilde{g}(e_a,X) + 2\d_{0a}(\tilde{g}(e_0,X) 
- \r^2Xf )\\ 
&=& 2\tilde{g}(e_a,X) + 2\d_{0a}(-2\r^2 \theta_0^P(X) - \r^2Xf) 
= 2\tilde{g}(e_a,X).
\end{eqnarray*}
Here we have used\re{ge0Equ} and $\theta_0^P=-\frac{1}{2}df$. 

Similarly, we get
\[ 2\tilde{g}(D_{X}e_0,e_a) = X\tilde{g}(e_0,e_a)  + e_0\tilde{g}(X,e_a) - 
e_a\tilde{g}(e_0,X) =  2 \tilde{g}(X,e_a).\] 
To finishes the proof of Theorem \ref{hKconifThm} it remains 
to compute the signature of $\hat{g}$. 
One can easily check that the matrix \re{tildegmatrixEqu}
has signature $(4,4)$ if $f_1h<0$, signature $(8,0)$ if
$h>0$ and  $f_1>0$ and signature $(0,8)$ if $h<0$ and $f_1<0$. 
This implies that $\hat{g}$ has signature $(4k,4\ell +4)$
if  $f_1<0$ and signature $(4k+4,4\ell )$ if $f_1>0$.  
\qed

Any conical pseudo-hyper-K\"ahler manifold $(M,g,J_1,J_2,J_3,\xi )$ 
is foliated by the leaves of the four-dimensional integrable distribution
defined by the vector fields $\xi, J_1\xi , J_2\xi , J_3\xi$.
The space of leaves inherits a quaternionic pseudo-K\"ahler
structure, at least if we restrict the foliation to 
a suitable open subset of $M$. Let us denote by $(\bar{M},\bar{g}, Q)$ 
the (at least locally defined) quaternionic pseudo-K\"ahler manifold 
associated 
with the conical pseudo-hyper-K\"ahler manifold $(\hat{M}, 
\hat{g},\hat{J}_1 ,\hat{J}_2 ,\hat{J}_3, \xi )$ of Theorem \ref{hKconifThm}.  

\bc \label{signatureCor} The signature of the 
quaternionic pseudo-K\"ahler manifold $(\bar{M},\bar{g}, Q)$ 
resulting from Theorem \ref{hKconifThm} 
depends only on the 
signature $(4k,4\ell )$ of the original pseudo-hyper-K\"ahler manifold 
$(M,g,J_1,J_2,J_3,Z)$ and on the signs of the functions $f$ 
and $f_1=f-h$, where $-f$ is the Hamiltonian chosen for the construction 
(unique up to an additive constant) and $h=g(Z,Z)/2$. It is 
$(4k,4\ell )$ if $f_1f>0$, 
$(4k-4,4\ell +4)$ if $f>0$ and $f_1<0$ and $(4k+4,4\ell -4)$ 
if $f<0$ and $f_1>0$. 
\ec

\pf This follows from the fact that the 
signature of $\bar{g}$ is obtained from that of $\hat{g}$ by subtracting
$(4,0)$ if $f>0$ and $(0,4)$ if $f<0$. 
\qed 
\bc \label{RiemCor} 
The construction of Theorem \ref{hKconifThm} yields a (positive 
definite) 
quaternionic K\"ahler manifold $(\bar{M},\bar{g}, Q)$ of
positive scalar curvature if and only if the metric $g$ of the 
original pseudo-hyper-K\"ahler manifold 
$(M,g,J_1,J_2,J_3,Z)$ is positive definite and $f_1>0$. 
It yields a quaternionic K\"ahler manifold $(\bar{M},\bar{g}, Q)$ of
negative scalar curvature if and only if $g$ 
is either positive definite and $f<0$ or 
if it has signature $(4k,4)$, $f<0$ and $f_1>0$. 
\ec 

\noindent 
{\bf Remark:} Notice that 
Theorem \ref{hKconifThm} provides us with a quaternionic 
pseudo-K\"ahler manifold for any choice of Hamiltonian $f$ for $Z$ 
with respect to $\o_1$ and any choice of 
$S^1$-principal bundle $(P,\eta )$ with 
connection $\eta$ such that the  curvature is $\o_1 -\frac{1}{2}d(g_N\zeta )$. 
Locally any two $S^1$-principal bundles with the same curvature
are equivalent such that, for the local geometry,  
the only essential choice is the Hamiltonian $f$, which is 
unique up to a constant $c$. It follows from 
\cite{HKLR} p.\ 553-554 that, up to a constant factor,
$f$ is a K\"ahler potential with respect to $J_2$.  
The Hamiltonian will be explicitly computed in 
the examples of the next section.  
\section{Application to the c-map}
\bd  \label{CASKDef}
A {\cmssl conical (affine) special K\"ahler manifold} $(M,J,g,\n ,\xi )$ is a 
pseudo-K\"ahler manifold $(M,J,g)$ endowed with a flat torsionfree 
connection $\n$ and a vector field $\xi$ such that 
\begin{enumerate}
\item[(i)] $\n \o=0$, where $\o=gJ$ is the K\"ahler form,
\item[(ii)] $d^\n J=0$, where $J$ is considered as a one-form with values
in the tangent bundle, 
\item[(iii)] $\n \xi = D\xi = {\rm Id}$, where $D$ is the Levi-Civita
connection,    
\item[(iv)] $g$ is definite on the plane $\mathcal{D} = \mathrm{span}
\{ \xi , J\xi \}$.   
\end{enumerate}
\ed
The above definition is slightly more general than the 
definition of a conical special K\"ahler manifold in \cite{CM}, for instance, 
since here we prefer not to restrict the signature of the metric.
The rigid c-map associates with $M$ the  
pseudo-hyper-K\"ahler manifold  $(N=T^*M,g_N,J_1,J_2,J_3)$, with 
the geometric data defined as follows, cf.\ \cite{ACD}. 
Using the connection $\n$ we can identify 
$TN = T^hN \oplus T^vN= \pi^*TM\oplus \pi^* T^*M$, 
where $\pi : N=T^*M\ra M$ is the canonical projection, 
$T^vN = \mathrm{ker}\, d\pi$ is the vertical 
distribution and 
$T^hN$ is the horizontal distribution defined by $\n$.   
Using these identifications, we have 
\[ g_N= \left( \begin{array}{cc}g&0\\
0&g^{-1}
\end{array}\right) , \quad J_1 = \left( \begin{array}{cc}J&0\\
0&J^*
\end{array}\right) ,\quad J_2 = \left( \begin{array}{cc}0&-\o^{-1}\\
\o &0
\end{array}\right) ,\quad J_3 = J_1J_2 .\]
The vector field $Z=J\xi$ is a Hamiltonian Killing vector field
on $M$. In fact, $DZ=JD\xi=J$ is skew-symmetric and the 
function $h=\frac{g(Z,Z)}{2}$ satisfies $dh= g(DZ,Z)=g(J\cdot , Z)=
-\o Z$. We extend $Z$ to a vector field $Z_N$ on $N$ by 
\[ Z_N(\pi^*q^i) := \pi^*Z(q^i),\quad Z_N(p_i) = 0, \]
where $(q^i)$ are $\n$-affine local coordinates on 
$M$ and $(\pi^*q^i,p_i)$ are the corresponding coordinates
of $N=T^*M$. One can easily check that this extension does not 
depend on the choice of affine coordinates.  
\bp For any pseudo-hyper-K\"ahler manifold  $(N,g_N,J_1,J_2,J_3)$ 
obtained from the rigid c-map, the 
vector field $Z_N$ is a time-like or space-like $\o_1$-Hamiltonian 
Killing vector field,  
which satisfies $D_{Z_N}Z_N=J_1Z_N$ and $\mathcal{L}_{Z_N}J_2=-J_3$. 
\ep 

\pf $Z$ is time-like or space-like by (iv) of Definition \ref{CASKDef}.
This implies that $Z_N$ is time-like or space-like. 
Let us recall 
that there exist $\n$-affine local coordinates 
$(q^i)$ on $M$ such that $\xi = \sum q^i\partial_i$, see \cite{CM1} Prop.\ 5. 
The following calculations will be always 
in such coordinates. Notice that $Z=J\xi = \sum J^i_jq^j\partial_i$. 
Since $\mathcal{L}_{Z}g=0$, we have 
\[ \mathcal{L}_{Z_N}g_N= \pi^* \mathcal{L}_{Z}g 
+  \sum Z(g^{ij})dp_idp_j= \sum Z(g^{ij})dp_idp_j.\] 
Since $Z(g^{ij})=-\sum g^{ik}Z(g_{kl})g^{lj}$, it suffices to show that 
$Z(g_{kl})=0$. Let us first recall\footnote{$(M,g,\n )$ is in fact an 
intrinsic affine hypersphere \cite{BC}, 
which implies the symmetry of $\n g$.}  
that $\n g$ is totally symmetric, 
which follows from Definition \ref{CASKDef} (i-ii) using the 
skew-symmetry of $J$. Using this property and Definition \ref{CASKDef} (i-ii),  
we obtain
\begin{eqnarray*} Z(g_{kl}) &=& \sum J^i_jq^jg_{kl,i} = \sum J^i_jq^jg_{ki,l} = 
-\sum q^j(J^i_j)_{,l} g_{ki}= -\sum q^j(J^i_l)_{,j} g_{ki}\\
&=& \sum q^jJ^i_l g_{ki,j}
 = \sum J^i_l\xi (g_{ki}).
\end{eqnarray*}
We claim that $\xi (g_{ki})=0$. 
Let us first observe that  
$D\xi =\mathrm{Id}$ 
implies $\mathcal{L}_\xi g = 2g$. This shows that $\xi (g_{ij})=0$, since
$\mathcal{L}_\xi q^i= q^i$. Summarizing, we have proven that $Z_N$
is a Killing vector field.
 
Next we prove that $Z_N$ is Hamiltonian with respect to $\o_1 =  g_NJ_1$. 
We consider the function $h=\frac{g(Z,Z)}{2}$. Then the 
calculation 
\be \label{HamEqu} d (\pi^*h) = \pi^*dh = -\pi^*(\o Z) = -\o_1 Z_N\ee 
proves that $-h$ is a Hamiltonian function for $Z_N$. This implies
the equation $D_{Z_N}Z_N=J_1Z_N$, as follows from Lemma \ref{dhLemma}.

Finally, we check that $\mathcal{L}_{Z_N} J_2 = - J_3$ or, equivalently,
that $\mathcal{L}_{Z_N}\o_2 = -\o_3$. Notice that
\[ J_3 = \left( \begin{array}{cc}0&-g^{-1}\\
g &0
\end{array}\right) .\] 
So 
\[ \o_2 =  \left( \begin{array}{cc}0& -J^*\\
J&0\end{array}\right) ,\quad 
\o_3 =  \left( \begin{array}{cc}0& -\id\\
\id &0\end{array}\right) ,\] 
that is   
\begin{eqnarray*} \o_2 &=& \sum J^j_idq^i\wedge dp_j \\
\o_3 &=& \sum dq^i\wedge dp_i .
\end{eqnarray*}
It is sufficient to check that $\mathcal{L}_{Z_N}\o_3 = \o_2$. 
Now 
\[ \mathcal{L}_{Z_N} \o_3 =  \sum d(Jq)^i\wedge dp_i = \o_2 + \sum 
J^i_{k,l}q^kdq^k\wedge dp_i =  \o_2, \] 
where $(Jq)^i= \sum J^i_jq^j$ and we have used that 
\[ \sum J^i_{k,l}q^k=
\sum q^kJ^i_{l,k}= \xi (J^i_l) = \xi (g^{ij} \o_{jl}) 
=0.\] 
\qed 
\bc For any pseudo-hyper-K\"ahler manifold  $(N,g_N,J_1,J_2,J_3)$ 
obtained from the rigid c-map, the assumptions for the 
conification construction of Theorem \re{hKconifThm} are satisfied 
for the Killing vector field $\zeta = 2Z_N$.  Therefore any choice of 
Hamiltonian $f$ for $\zeta$ with respect to $\o_1$ and any choice 
of $S^1$-principal bundle $(P,\eta )$ with connection $\eta$ 
with the curvature $\o_1 -\frac{1}{2}d(g_N\zeta )$ defines a 
conical  pseudo-hyper-K\"ahler manifold $(\hat{N},g_{\hat{N}},\hat{J}_1,
\hat{J}_2, \hat{J}_3, \xi )$.  
\ec 

Consider now the function $h := \frac{1}{2}g_N(\zeta , \zeta )
= 2g_N(Z_N,Z_N)$ 
associated with the rescaled vector field $\zeta = 2Z_N$. Then 
\re{HamEqu} shows that 
any function $f$ satisfying $df=-\o_1 \zeta$ is of the form 
\[ f= \frac{1}{2}h + c,\]
for some constant $c$. The choice $c=0$ will be called the 
\emph{canonical choice of Hamiltonian}. 
Thus the function 
$f_1 := f-h$ is now given by 
\[ f_1= -\frac{1}{2}h +c.\]

Now we specialise to the cases where the resulting quaternionic K\"ahler 
manifold $(\bar{N},\bar{g},Q)$ is positive definite, see 
Corollary \ref{RiemCor}.  
\bc If the conical special K\"ahler manifold $M$ is positive definite, 
then the resulting quaternionic 
K\"ahler manifold $\bar{N}$ has positive scalar curvature if 
$-\frac{1}{2}h +c>0$ and negative scalar curvature if $\frac{1}{2}h +c<0$.
If the conical special K\"ahler manifold $M$ has signature $(2k,2)$
with time-like Euler field and if  $-\frac{1}{2}h +c>0$ and 
$\frac{1}{2}h +c<0$, then 
$\bar{N}$ has negative scalar curvature. 
In particular, for the 
canonical choice of Hamiltonian the scalar curvature of the quaternionic 
K\"ahler manifold $\bar{N}$ is always negative. 
\ec 
The last result is consistent with our conjecture that 
the canonical choice of Hamiltonian yields the Ferrara-Sabharwal
metric (up to a factor). The deformation by the constant $c$ should correspond
to the one-loop correction of the metric considered in \cite{APSV}. 
The determination of the precise relation between the constant $c$ and the
one-loop parameter is left for the future.

\end{document}